\documentclass[12pt]{article}

\usepackage{amsmath,amssymb,amsthm,amscd,a4wide,upref}

\usepackage[enableskew]{youngtab}
\usepackage{ytableau}

\usepackage{hyperref}
\hypersetup{colorlinks,citecolor=blue,filecolor=black,linkcolor=blue,urlcolor=blue}

\usepackage{lmodern}     
\usepackage[T1]{fontenc}

\begin{document}


\newcommand{\non}{\nonumber}
\newcommand{\scl}{\scriptstyle}
\newcommand{\sclnearrow}{{\scl\nearrow}\ts}
\newcommand{\scloplus}{{\scl\bigoplus}}
\newcommand{\wt}{\widetilde}
\newcommand{\wh}{\widehat}
\newcommand{\ot}{\otimes}
\newcommand{\fand}{\quad\text{and}\quad}
\newcommand{\Fand}{\qquad\text{and}\qquad}
\newcommand{\ts}{\,}
\newcommand{\tss}{\hspace{1pt}}
\newcommand{\lan}{\langle\ts}
\newcommand{\ran}{\ts\rangle}
\newcommand{\vl}{\tss|\tss}
\newcommand{\qin}{q^{-1}}
\newcommand{\tpr}{t^{\tss\prime}}
\newcommand{\spr}{s^{\tss\prime}}
\newcommand{\di}{\partial}
\newcommand{\hra}{\hookrightarrow}
\newcommand{\antiddots}
    {\underset{\displaystyle\cdot\quad\ }
    {\overset{\displaystyle\quad\ \cdot}{\cdot}}}
\newcommand{\dddots}
    {\underset{\displaystyle\quad\ \cdot}
    {\overset{\displaystyle\cdot\quad\ }{\cdot}}}
\newcommand{\atopn}[2]{\genfrac{}{}{0pt}{}{#1}{#2}}

\newcommand{\su}{s^{}}
\newcommand{\vac}{\mathbf{1}}
\newcommand{\vacf}{\tss|0\rangle}
\newcommand{\BL}{ {\overline L}}
\newcommand{\BE}{ {\overline E}}
\newcommand{\BP}{ {\overline P}}
\newcommand{\ol}{\overline}
\newcommand{\pr}{^{\tss\prime}}
\newcommand{\ba}{\bar{a}}
\newcommand{\bb}{\bar{b}}
\newcommand{\eb}{\bar{e}}
\newcommand{\bi}{\bar{\imath}}
\newcommand{\bj}{\bar{\jmath}}
\newcommand{\bk}{\bar{k}}
\newcommand{\bl}{\bar{l}}
\newcommand{\hb}{\mathbf{h}}
\newcommand{\gb}{\mathbf{g}}
\newcommand{\For}{\qquad\text{or}\qquad}
\newcommand{\OR}{\qquad\text{or}\qquad}
\newcommand{\emp}{\mbox{}}


\newcommand{\U}{{\rm U}}
\newcommand{\Z}{{\rm Z}}
\newcommand{\ZY}{{\rm ZY}}
\newcommand{\Ar}{{\rm A}}
\newcommand{\Br}{{\rm B}}
\newcommand{\Cr}{{\rm C}}
\newcommand{\Fr}{{\rm F}}
\newcommand{\Mr}{{\rm M}}
\newcommand{\Sr}{{\rm S}}
\newcommand{\Prm}{{\rm P}}
\newcommand{\Lr}{{\rm L}}
\newcommand{\oL}{{\ol L}}
\newcommand{\oM}{{\ol M}}
\newcommand{\Ir}{{\rm I}}
\newcommand{\Jr}{{\rm J}}
\newcommand{\Qr}{{\rm Q}}
\newcommand{\Rr}{{\rm R}}
\newcommand{\X}{{\rm X}}
\newcommand{\Y}{{\rm Y}}
\newcommand{\DY}{ {\rm DY}}
\newcommand{\Or}{{\rm O}}
\newcommand{\SO}{{\rm SO}}
\newcommand{\GL}{{\rm GL}}
\newcommand{\Spr}{{\rm Sp}}
\newcommand{\Zr}{{\rm Z}}
\newcommand{\ev}{{\rm ev}}
\newcommand{\Pf}{{\rm Pf}}
\newcommand{\Ann}{{\rm{Ann}\ts}}
\newcommand{\Norm}{{\rm Norm\tss}}
\newcommand{\Ad}{{\rm Ad}}
\newcommand{\SY}{{\rm SY}}
\newcommand{\Pff}{{\rm Pf}\tss}
\newcommand{\Hf}{{\rm Hf}\tss}
\newcommand{\trts}{{\rm tr}\ts}
\newcommand{\otr}{{\rm otr}}
\newcommand{\row}{{\rm row}}
\newcommand{\End}{{\rm{End}\ts}}
\newcommand{\Mat}{{\rm{Mat}}}
\newcommand{\Hom}{{\rm{Hom}}}
\newcommand{\id}{{\rm id}}
\newcommand{\middd}{{\rm mid}}
\newcommand{\ch}{{\rm{ch}\ts}}
\newcommand{\ind}{{\rm{ind}\ts}}
\newcommand{\Normts}{{\rm{Norm}\ts}}
\newcommand{\mult}{{\rm{mult}}}
\newcommand{\per}{{\rm per}\ts}
\newcommand{\sgn}{{\rm sgn}\ts}
\newcommand{\sign}{{\rm sign}\ts}
\newcommand{\qdet}{{\rm qdet}\ts}
\newcommand{\sdet}{{\rm sdet}\ts}
\newcommand{\Ber}{{\rm Ber}\ts}
\newcommand{\inv}{{\rm inv}\ts}
\newcommand{\inva}{{\rm inv}}
\newcommand{\grts}{{\rm gr}\ts}
\newcommand{\grpr}{{\rm gr}^{\tss\prime}\ts}
\newcommand{\degpr}{{\rm deg}^{\tss\prime}\tss}
\newcommand{\Cond}{ {\rm Cond}\tss}
\newcommand{\Fun}{{\rm{Fun}\ts}}
\newcommand{\Rep}{{\rm{Rep}\ts}}
\newcommand{\sh}{{\rm{sh}}}
\newcommand{\weight}{{\rm{wt}\ts}}
\newcommand{\chara}{{\rm{char}\ts}}
\newcommand{\diag}{ {\rm diag}}
\newcommand{\Bos}{ {\rm Bos}}
\newcommand{\Ferm}{ {\rm Ferm}}
\newcommand{\cdet}{ {\rm cdet}}
\newcommand{\rdet}{ {\rm rdet}}
\newcommand{\imm}{ {\rm imm}}
\newcommand{\ad}{ {\rm ad}}
\newcommand{\tr}{ {\rm tr}}
\newcommand{\gr}{ {\rm gr}\tss}
\newcommand{\str}{{\rm str^{}}}
\newcommand{\loc}{{\rm loc}}
\newcommand{\Gr}{{\rm G}}

\newcommand{\twobar}{{\bar 2}}
\newcommand{\threebar}{{\bar 3}}


\newcommand{\AAb}{\mathbb{A}\tss}
\newcommand{\CC}{\mathbb{C}}
\newcommand{\KK}{\mathbb{K}\tss}
\newcommand{\QQ}{\mathbb{Q}\tss}
\newcommand{\SSb}{\mathbb{S}\tss}
\newcommand{\TT}{\mathbb{T}\tss}
\newcommand{\ZZ}{\mathbb{Z}\tss}
\newcommand{\Sbb}{\mathbb{S}}
\newcommand{\ZZb}{\mathbb{Z}}


\newcommand{\Ac}{{\mathcal A}}
\newcommand{\Bc}{{\mathcal B}}
\newcommand{\Cc}{{\mathcal C}}
\newcommand{\Cl}{{\mathcal Cl}}
\newcommand{\Dc}{{\mathcal D}}
\newcommand{\Ec}{{\mathcal E}}
\newcommand{\Fc}{{\mathcal F}}
\newcommand{\Jc}{{\mathcal J}}
\newcommand{\Gc}{{\mathcal G}}
\newcommand{\Hc}{{\mathcal H}}
\newcommand{\Lc}{{\mathcal L}}
\newcommand{\Nc}{{\mathcal N}}
\newcommand{\Xc}{{\mathcal X}}
\newcommand{\Yc}{{\mathcal Y}}
\newcommand{\Oc}{{\mathcal O}}
\newcommand{\Pc}{{\mathcal P}}
\newcommand{\Qc}{{\mathcal Q}}
\newcommand{\Rc}{{\mathcal R}}
\newcommand{\Sc}{{\mathcal S}}
\newcommand{\cS}{{\check S}}
\newcommand{\Tc}{{\mathcal T}}
\newcommand{\Uc}{{\mathcal U}}
\newcommand{\Vc}{{\mathcal V}}
\newcommand{\Wc}{{\mathcal W}}
\newcommand{\Zc}{{\mathcal Z}}
\newcommand{\HC}{{\mathcal {HC}}}


\newcommand{\asf}{\mathsf a}
\newcommand{\bsf}{\mathsf b}
\newcommand{\csf}{\mathsf c}
\newcommand{\nsf}{\mathsf n}


\newcommand{\Sym}{\mathfrak S}
\newcommand{\h}{\mathfrak h}
\newcommand{\q}{\mathfrak q}
\newcommand{\n}{\mathfrak n}
\newcommand{\m}{\mathfrak m}
\newcommand{\p}{\mathfrak p}
\newcommand{\gl}{\mathfrak{gl}}
\newcommand{\oa}{\mathfrak{o}}
\newcommand{\spa}{\mathfrak{sp}}
\newcommand{\osp}{\mathfrak{osp}}
\newcommand{\g}{\mathfrak{g}}
\newcommand{\kgot}{\mathfrak{k}}
\newcommand{\agot}{\mathfrak{a}}
\newcommand{\bgot}{\mathfrak{b}}
\newcommand{\sll}{\mathfrak{sl}}
\newcommand{\f}{\mathfrak{f}}
\newcommand{\z}{\mathfrak{z}}
\newcommand{\Zgot}{\mathfrak{Z}}


\newcommand{\al}{\alpha}
\newcommand{\be}{\beta}
\newcommand{\ga}{\gamma}
\newcommand{\de}{\delta}
\newcommand{\De}{\Delta}
\newcommand{\Ga}{\Gamma}
\newcommand{\ep}{\epsilon}
\newcommand{\ee}{\epsilon^{}}
\newcommand{\ve}{\varepsilon}
\newcommand{\ls}{\ts\lambda\ts}
\newcommand{\vk}{\varkappa}
\newcommand{\vs}{\varsigma}
\newcommand{\vt}{\vartheta}
\newcommand{\ka}{\kappa}
\newcommand{\vp}{\varphi}
\newcommand{\la}{\lambda}
\newcommand{\La}{\Lambda}
\newcommand{\si}{\sigma}
\newcommand{\ze}{\zeta}
\newcommand{\om}{\omega}
\newcommand{\Om}{\Omega}
\newcommand{\up}{\upsilon}


\newtheorem{thm}{Theorem}[section]
\newtheorem{lemma}[thm]{Lemma}
\newtheorem{prop}[thm]{Proposition}
\newtheorem{cor}[thm]{Corollary}
\newtheorem{conj}[thm]{Conjecture}

\theoremstyle{definition}
\newtheorem{definition}[thm]{Definition}
\newtheorem{example}[thm]{Example}

\theoremstyle{remark}
\newtheorem{remark}[thm]{Remark}

\newcommand{\bth}{\begin{thm}}
\renewcommand{\eth}{\end{thm}}
\newcommand{\bpr}{\begin{prop}}
\newcommand{\epr}{\end{prop}}
\newcommand{\ble}{\begin{lemma}}
\newcommand{\ele}{\end{lemma}}
\newcommand{\bco}{\begin{cor}}
\newcommand{\eco}{\end{cor}}
\newcommand{\bex}{\begin{example}}
\newcommand{\eex}{\end{example}}
\newcommand{\bde}{\begin{definition}}
\newcommand{\ede}{\end{definition}}
\newcommand{\bre}{\begin{remark}}
\newcommand{\ere}{\end{remark}}
\newcommand{\bcj}{\begin{conj}}
\newcommand{\ecj}{\end{conj}}

\renewcommand{\theequation}{\arabic{section}.\arabic{equation}}

\numberwithin{equation}{section}


\newcommand{\bpf}{\begin{proof}}
\newcommand{\epf}{\end{proof}}


\def\beql#1{\begin{equation}\label{#1}}

\newcommand{\bal}{\begin{aligned}}
\newcommand{\eal}{\end{aligned}}
\newcommand{\beq}{\begin{equation}}
\newcommand{\eeq}{\end{equation}}
\newcommand{\ben}{\begin{equation*}}
\newcommand{\een}{\end{equation*}}

\title{\Large\bf Central elements and evaluation map for\\
the quantum queer superalgebras}

\author{{Ming Liu,\quad Alexander Molev\quad and\quad Jian Zhang}}

\date{} 
\maketitle


\begin{abstract}
We consider the $R$-matrix presentations of the quantum queer superalgebra $\U_q(\q_n)$
and its affine counterpart $\U_q(\wh\q_n)$. We derive crossing symmetry relations
for the $R$-matrices and use them to
construct central elements in both superalgebras. We also produce an epimorphism
$\ev:\U_q(\wh\q_n)\to\U_q(\q_n)$ identical on the subalgebra
isomorphic to $\U_q(\q_n)$.



\end{abstract}

\section{Introduction}
\label{sec:int}

A quantum deformation $\U_q(\q_n)$ of the universal enveloping algebra for the queer Lie
superalgebra $\q_n$ was introduced by Olshanski~\cite{o:qu}. He described its algebraic structure
by proving the Poincar\'{e}--Birkhoff--Witt theorem and developed a quantum analogue of the Schur--Sergeev
duality between $\U_q(\q_n)$ and the Hecke--Clifford superalgebra which was also introduced in \cite{o:qu}.

A Drinfeld-type presentation of $\U_q(\q_n)$ was produced by Grantcharov {\em et al.} in \cite{gjkk:hw},
where classification results on the highest weight modules were also obtained.
For a further development of the representation theory of $\U_q(\q_n)$,
including the crystal basis theory, see \cite{gjkkk:cb}.

An affine counterpart of $\U_q(\q_n)$ was introduced by Chen and Guay~\cite{cg:ta}
under the name {\em twisted affine quantum superalgebra}.
The name is justified by the fact that it is
a quantization of a Lie bialgebra structure on the twisted affine superalgebra corresponding
to $\q_n$. We will only use its version with the trivial central charge
and call it
the {\em quantum queer loop superalgebra} to be denoted by $\U_q(\wh\q_n)$. Since
no non-twisted version of this superalgebra is known, no confusion
should occur; this is also consistent with the Yangian terminology used in \cite{n:yq}.
One of the main results of \cite{cg:ta} is the construction of
a functor of Schur--Weyl-type connecting $\U_q(\wh\q_n)$ to the affine Hecke--Clifford superalgebras
of Jones and Nazarov~\cite{jn:as}, and the proof, that
it provides an equivalence between two categories of modules.

Our first main result is a construction of the epimorphism
\beql{ev}
\ev:\U_q(\wh\q_n)\to\U_q(\q_n).
\eeq
The quantum queer superalgebra $\U_q(\q_n)$ turns out to be
a natural subalgebra of $\U_q(\wh\q_n)$ and the epimorphism \eqref{ev}
is identical on this subalgebra. The epimorphism allows
one to construct evaluation modules $V$ over $\U_q(\wh\q_n)$ by extending the action
of $\U_q(\q_n)$ on $V$. The use of the Hopf algebra
structure on $\U_q(\wh\q_n)$ then produces a wide class of $\U_q(\wh\q_n)$-modules.

Furthermore, we construct families of central elements in both superalgebras
$\U_q(\q_n)$ and $\U_q(\wh\q_n)$. In the first case,
the central elements of one of the families coincide with those
already
constructed in \cite{s:qi}
via a categorical approach; see \eqref{savagesrtce}
below.\footnote{We thank Alistair Savage for pointing out the connection.}
As suggested in \cite[Sec.~8]{s:qi},
it is likely that together with one additional element they generate the center
of $\U_q(\q_n)$.

In the affine case, the central elements occur as the coefficients
of a series $z(u)$, whose analogues for the quantum affine algebra $\U_q(\wh\gl_n)$
were found in \cite{br:nb}; see also \cite{jlm:eq}.
It is reasonable to expect that the families of central
elements in $\U_q(\q_n)$ and $\U_q(\wh\q_n)$ are related by
the epimorphism \eqref{ev}; see Example~\ref{ex:none}.

\section{Evaluation homomorphism}
\label{sec:eh}

Consider the $\ZZ_2$-graded vector space $\CC^{n|n}$ over the field of complex numbers with the
canonical basis
$e_{-n},\dots,e_{-1},e_1,\dots,e_n$, where
the vector $e_i$ has the parity
$\bi\mod 2$ and
\ben
\bi=\begin{cases} 1\qquad\text{for}\quad i<0,\\
0\qquad\text{for}\quad i>0.
\end{cases}
\een
The endomorphism algebra $\End\CC^{n|n}$ is then equipped with a $\ZZ_2$-gradation with
the parity of the matrix unit $E_{ij}$ found by
$\bi+\bj\mod 2$. We will identify
the algebra of
even matrices over a superalgebra $\Ac$ with the tensor product algebra
$\Ac\ot\End\CC^{n|n}$, so that a square matrix $A=[A_{ij}]$ of size $2n$
is regarded as the element
\beql{signsma}
A=\sum_{i,j=-n}^{n}A_{ij}\ot E_{ij}(-1)^{\bi\tss\bj+\bj}\in\Ac\ot\End\CC^{n|n},
\eeq
where the entries $A_{ij}$ are assumed to be homogeneous of parity $\bi+\bj\mod 2$.
The extra signs are necessary to keep the usual rule for the matrix multiplication
without signs.
The zero value will be skipped in similar summation formulas.

The {\em Olshanski $R$-matrix} is the element $S\in \End\CC^{n|n}\ot\End\CC^{n|n}$
defined by
\ben
S=\sum_{i,j=-n}^n q^{(\de_{i,j}+\de_{-i,j})(1-2\bj)}_{}E_{ii}\ot E_{jj}
+\ve \sum_{i>j}(-1)^{\bj}(E_{ij}+E_{-i,-j})\ot E_{ji},
\een
where
$q\in\CC$ is nonzero and $\ve=q-\qin$. As shown in \cite{o:qu}, the element $S$
satisfies the Yang--Baxter equation.

Following \cite{o:qu},
define the {\em quantum queer superalgebra} $\U_q(\q_n)$
as a $\ZZ_2$-graded algebra generated by the entries $L_{ij}$
of the upper-triangular matrix $L=[L_{ij}]$.
The generators
$L_{ij}$ are of the parity $\bi+\bj$, where $i\leqslant j$ with $i$ and $j$
running over the set $\{-n,\dots,-1,1,\dots,n\}$.
They
satisfy the following defining relations.
Consider the elements of the tensor product superalgebra
$\U_q(\q_n)\ot\End\CC^{n|n}\ot\End\CC^{n|n}$ given by
\ben
L_1=\sum_{i\leqslant j} L_{ij}\ot E_{ij}\ot 1\fand
L_2=\sum_{i\leqslant j} L_{ij}\ot 1\ot E_{ij},
\een
with no signs occurring in accordance with \eqref{signsma}.
The defining relations for the superalgebra $\U_q(\q_n)$ take
the form of the relation
\beql{RTT}
L_1 L_2\tss S=S\tss L_2 L_1,
\eeq
where $S$ is identified with the element $1\ot S$, together with
the relations
\ben
L_{ii}L_{-i,-i}=L_{-i,-i}L_{ii}=1,\qquad i=1,\dots,n.
\een

Now consider the {\em Jones--Nazarov $R$-matrix}
\beql{jns}
S(u,v)=S+\frac{\ve\tss P}{u^{-1}v-1}+\frac{\ve\tss PJ_1J_2}{u\tss v-1},
\eeq
where the permutation operator $P$ is given by
\ben
P=\sum_{i,j=-n}^n E_{ij}\ot E_{ji}(-1)^{\bj}\in \End\CC^{n|n}\ot\End\CC^{n|n},
\een
while $J_1=J\ot 1$ and $J_2=1\ot J$ for the odd element $J$ given by
\ben
J=\sum_{i=-n}^n E_{i,-i}(-1)^{\bi}\in \End\CC^{n|n}.
\een
According to
\cite{jn:as}, the function $S(u,v)$
satisfies the parameter-dependent Yang--Baxter equation.

Following Chen and Guay~\cite{cg:ta},
define the {\em quantum queer loop superalgebra} $\U_q(\wh\q_n)$
as a $\ZZ_2$-graded algebra with generators
$L^{(r)}_{ij}$ of parity $\bi+\bj$, where $i$ and $j$
run over the set of indices $\{-n,\dots,-1,1,\dots,n\}$ and $r=0,1,\dots$.
The generators
satisfy the following defining relations.
Introduce the formal series
\beql{tiju}
L_{ij}(u)=\sum_{r=0}^{\infty}L_{ij}^{(r)}\ts u^{-r}
\in \U_q(\wh\q_n)[[u^{-1}]]
\eeq
and combine them into the matrix $L(u)=[L_{ij}(u)]$.
Consider the elements of the tensor product superalgebra
$\U_q(\wh\q_n)\ot\End\CC^{n|n}\ot\End\CC^{n|n}$ given by
\ben
L_1(u)=\sum_{i,j=-n}^{n} L_{ij}(u)\ot E_{ij}\ot 1(-1)^{\bi\tss\bj+\bj}\fand
L_2(u)=\sum_{i,j=-n}^{n} L_{ij}(u)\ot 1\ot E_{ij}(-1)^{\bi\tss\bj+\bj}.
\een
The defining relations for the superalgebra $\U_q(\wh\q_n)$ take
the form of the relation
\beql{RTTuv}
L_1(u)\tss L_2(v)\tss S(u,v)=S(u,v)\tss L_2(v)\tss L_1(u),
\eeq
together with the conditions on the zero mode matrix $L^{(0)}=[L^{(0)}_{ij}]$:
\ben
L^{(0)}_{ii}L^{(0)}_{-i,-i}=L^{(0)}_{-i,-i}L^{(0)}_{ii}=1,\qquad i=1,\dots,n,
\een
and
\ben
L^{(0)}_{ij}=0\qquad\text{for}\quad i>j.
\een

It is immediate from the defining relations of both superalgebras that the mapping
\beql{homim}
L\to L^{(0)}
\eeq
defines a homomorphism $\U_q(\q_n)\to \U_q(\wh\q_n)$.

To state our first main result, introduce the even matrix $\oL$ by
\beql{lbar}
\oL=JLJ=-\sum_{i\leqslant j} L_{ij}\ot E_{-i,-j}=-\sum_{i\geqslant j} L_{-i,-j}\ot E_{ij}.
\eeq

\bth\label{thm:eval}
The mapping
\beql{evexp}
\ev:L(u)\mapsto L+\oL\tss u^{-1}
\eeq
defines a homomorphism $\U_q(\wh\q_n)\to \U_q(\q_n)$.
\eth

\bpf
We need to verity that the relation
\beql{imRTT}
(L_1+\oL_1\tss u^{-1})\tss (L_2+\oL_2\tss v^{-1})\tss S(u,v)
=S(u,v)\tss (L_2+\oL_2\tss v^{-1})\tss (L_1+\oL_1\tss u^{-1})
\eeq
holds in $\U_q(\q_n)\ot\End\CC^{n|n}\ot\End\CC^{n|n}$.

\ble\label{lem:relq}
We have the relations in $\U_q(\q_n)\ot\End\CC^{n|n}\ot\End\CC^{n|n}${\rm :}
\begin{align}
\oL_1 L_2\tss S&=S\tss L_2 \oL_1,
\label{olls}\\
L_1 \oL_2\tss (S-\ve P)&=(S-\ve P)\tss \oL_2L_1,
\label{lolsp}\\
\oL_1 \oL_2\tss (S-\ve PJ_1J_2)&=(S-\ve PJ_1J_2)\tss \oL_2\oL_1,
\label{ololsp}
\end{align}
together with
\begin{align}
\oL_1 L_2\tss PJ_1J_2&=PJ_1J_2\tss L_2\oL_1,\qquad L_1 \oL_2\tss PJ_1J_2=PJ_1J_2\tss \oL_2L_1,
\label{ollp}\\
L_1 L_2\tss PJ_1J_2&=PJ_1J_2\tss \oL_2\oL_1,\qquad \oL_1 \oL_2\tss PJ_1J_2=PJ_1J_2\tss L_2L_1.
\label{llpjj}
\end{align}
\ele

\bpf
We will make use of the
homomorphism from the Hecke--Clifford superalgebra $\HC_q(N)$,
\beql{homhc}
\HC_q(N)\to (\End\CC^{n|n})^{\ot N}
\eeq
as constructed in \cite[Theorem~5.2]{o:qu}. We only need the particular case $N=2$,
where $\HC_q(2)$ is generated by elements $T, c_1,c_2$ subject to the relations
\ben
T^2=\ve T+1,\qquad Tc_1=c_2T,\qquad Tc_2=c_1T-\ve(c_1-c_2),
\een
and
\ben
c_1^2=c_2^2=-1,\qquad c_1c_2+c_2c_1=0.
\een
The homomorphism \eqref{homhc} with $N=2$ takes the form
\beql{homnt}
T\mapsto PS,\qquad c_1\mapsto J_1,\qquad c_2\mapsto J_2.
\eeq
Hence, we derive the relations
\beql{js}
J_1S=SJ_1,\qquad J_2(S-\ve P)=(S-\ve P)J_2.
\eeq
Now \eqref{olls} follows from \eqref{RTT} by writing
\begin{multline}
\oL_1 L_2\tss S=J_1L_1J_1L_2\tss S=J_1L_1L_2J_1\tss S=J_1L_1L_2\tss SJ_1\\
{}=J_1SL_2L_1J_1=S J_1L_2L_1J_1=SL_2 J_1L_1J_1=S\tss L_2 \oL_1.
\non
\end{multline}
To verify \eqref{lolsp}, use the second relation in \eqref{js} to
perform a similar calculation by noting that
\ben
L_1 L_2\tss (S-\ve P)=(S-\ve P)\tss L_2 L_1.
\een
Relation \eqref{ololsp} is verified in the same way with the use
of the first relation in \eqref{llpjj}. In their turn, all
cases of \eqref{ollp} and \eqref{llpjj}
are easily checked.
\epf

After multiplying both sides of \eqref{imRTT}
by the denominators in \eqref{jns} and suitable powers
of $u$ and $v$, this becomes a relation for polynomials in $u$ and $v$.
It is straightforward to verify that all coefficients on both sides agree
due to \eqref{RTT} and additional relations recorded in Lemma~\ref{lem:relq}.
\epf

\bco\label{cor:emb}
The homomorphism \eqref{homim} defines an embedding
$\U_q(\q_n)\hra\U_q(\wh\q_n)$.
\eco

\bpf
The homomorphism \eqref{homim} is injective since its composition
with the evaluation homomorphism \eqref{evexp} is the identity map.
\epf

Corollary~\ref{cor:emb} allows us to regard $\U_q(\q_n)$
as a subalgebra of $\U_q(\wh\q_n)$.

Note that by twisting the homomorphism \eqref{evexp} with the sign
automorphism of  $\U_q(\wh\q_n)$ defined by $L(u)\mapsto L(-u)$,
we get another version of the evaluation homomorphism:
\ben
\ev':L(u)\mapsto L-\oL\tss u^{-1}.
\een

\section{Reflection equation superalgebra}
\label{sec:re}

Keeping notation \eqref{lbar}, introduce the matrix $M=L\oL^{-1}$.

\bpr\label{prop:refl}
The matrix $M$ satisfies the symmetry relation
\beql{symrel}
MJMJ=-1
\eeq
in $\U_q(\q_n)\ot\End\CC^{n|n}$,
and the
reflection equation
\beql{refleq}
M_1SM_2(S_{21}+\ve P J_1J_2)=SM_2(S_{21}+\ve P J_1J_2)M_1
\eeq
in $\U_q(\q_n)\ot\End\CC^{n|n}\ot\End\CC^{n|n}$,
where $S_{21}=PSP$. Equivalently, relation \eqref{refleq} can be written as
\beql{eqfore}
M_2\cS M_2(\cS-\ve J_1J_2)=\cS M_2(\cS-\ve J_1J_2)M_2,
\eeq
where $\cS=PS$.
\epr

\bpf
Relation \eqref{symrel} is immediate from the definition of the matrix $M$.
Clearly, \eqref{eqfore} follows from \eqref{refleq} by multiplying
both sides by $P$ from the right and from the left. To prove \eqref{refleq},
observe that
\beql{spmin}
(S-\ve P-\ve PJ_1J_2)(S_{21}+\ve P J_1J_2)=1.
\eeq
Indeed, using the homomorphism \eqref{homnt} we find that
\beql{ssto}
SS_{21}=\ve SP+1.
\eeq
Together with \eqref{js} this implies \eqref{spmin}. Now write
the left hand side of \eqref{refleq} as
\ben
L_1\oL_1^{-1}SL_2\oL_2^{-1}(S-\ve P-\ve PJ_1J_2)^{-1}
\een
which equals
\beql{llslim}
L_1L_2S\oL_1^{-1}\oL_2^{-1}(S-\ve P-\ve PJ_1J_2)^{-1}
\eeq
by \eqref{olls}. Furthermore, by \eqref{ololsp} we have
\ben
\oL_1 \oL_2\tss (S-\ve P-\ve PJ_1J_2)=(S-\ve P-\ve PJ_1J_2)\tss \oL_2\oL_1,
\een
hence, by inverting both sides and using \eqref{RTT}
we bring \eqref{llslim} to the form
\beql{sllsp}
SL_2L_1(S-\ve P-\ve PJ_1J_2)^{-1}\oL_2^{-1}\oL_1^{-1}.
\eeq
Finally, by \eqref{olls} and the first relation in \eqref{ollp}, we have
\ben
\oL_1L_2(S-\ve PJ_1J_2)=(S-\ve PJ_1J_2)L_2\oL_1.
\een
By conjugating both sides by $P$ we derive
\ben
\oL_2L_1(S_{21}+\ve PJ_1J_2)=(S_{21}+\ve PJ_1J_2)L_1\oL_2.
\een
Therefore, applying \eqref{spmin}, we can write \eqref{sllsp} as
\ben
SL_2L_1(S_{21}+\ve PJ_1J_2)\oL_2^{-1}\oL_1^{-1}=SL_2\oL_2^{-1}(S_{21}+\ve PJ_1J_2)L_1\oL_1^{-1}
\een
which coincides with the right hand side of \eqref{refleq}.
\epf

We will denote by $\U^{\circ}_q(\q_n)$ the subalgebra
of $\U_q(\q_n)$ generated by the entries of the matrix $M$.
We will also define an abstract superalgebra, regarding
\eqref{refleq} as the set of defining relations.

\bde\label{def:res}
The {\em reflection equation superalgebra} $\U_q^{\rm RE}(\q_n)$ associated with $\q_n$ is generated
by the entries $M_{ij}$ of the (abstract) matrix $M=[M_{ij}]$ subject to
the defining relations \eqref{symrel} and \eqref{refleq}.
\qed
\ede

By Proposition~\ref{prop:refl}, we have the homomorphism $\U_q^{\rm RE}(\q_n)\to \U^{\circ}_q(\q_n)$
taking the generator $M_{ij}$ to the $(i,j)$ entry of the matrix $L\oL^{-1}$.

\section{Central elements in $\U_q(\q_n)$}
\label{sec:cent-fin}

Introduce the diagonal matrix $D$ by
\beql{d}
D=\sum_{i=-n}^n q^{-2\tss|i|}\tss E_{ii}\in \End\CC^{n|n},
\eeq
where, as before, the zero value is skipped in the summation. We will keep using the matrix
$M=L\oL^{-1}$ introduced in the previous section. Define the {\em supertrace}
as the linear map
\ben
\str:\End\CC^{n|n}\to\CC,\qquad E_{ij}\mapsto \de_{ij}(-1)^{\bi}.
\een
The central elements provided by the following theorem
were obtained earlier in \cite{s:qi} via
a categorical approach in a slightly different form;
see Remark~\ref{rem:secre}(i) below. They
are analogous to
those found in \cite{rtf:ql} for the quantized enveloping algebra
$\U_q(\gl_n)$.

\bth\label{thm:centfin}
All elements
\beql{srtce}
\str\ts DM^k,\qquad k=1,2,\dots
\eeq
belong to the center of the superalgebra $\U_q(\q_n)$.
\eth

\bpf
In the tensor product superalgebra
$\U_q(\q_n)\ot\End\CC^{n|n}\ot\End\CC^{n|n}$
we have
\ben
L_2M_1=L_2L_1\oL_1^{-1}=S^{-1}L_1L_2S\oL_1^{-1}=S^{-1}L_1\oL_1^{-1}SL_2=S^{-1}M_1SL_2,
\een
where we used the consequence $L_2S\oL_1^{-1}=\oL_1^{-1}SL_2$ of \eqref{olls}.
An easy induction yields
\ben
L_2D_1M_1^k=D_1 L_2M_1^k=D_1S^{-1}M_1^kSL_2
\een
for any $k=1,2,\dots$. Hence the desired property of the elements \eqref{srtce}
will follow from the relation
\beql{sutre}
\str_1 D_1S^{-1}M_1^kS=\str_1 D_1M_1^k,
\eeq
where the subscript $1$ means that the supertrace is taken over
the first copy of the endomorphism superalgebra.

Consider the {\em super-transposition} $t$ which is defined as the antiautomorphism
of the superalgebra $\End\CC^{n|n}$ given by
\ben
t:E_{ij}\mapsto E_{ji}(-1)^{\bi\bj+\bj}.
\een
We will use a subscript as in $t_1,t_2$, etc. to indicate the copy of
the endomorphism superalgebra in multiple tensor products where the
super-transposition is applied.
Transform the left hand side of \eqref{sutre} by applying
the super-transposition $t_1$ to get
\ben
\str_1 D_1S^{-1}M_1^kS=\str_1 (D_1S^{-1})^{t_1}(M_1^kS)^{t_1}=
\str_1 (S^{-1})^{t_1}D_1S^{t_1}(M_1^k)^{t_1}.
\een
Since
\ben
\str_1 D_1(M_1^k)^{t_1}=\str_1 D_1M_1^k,
\een
the proof of relation
\eqref{sutre} will be completed by using
the following {\em crossing symmetry property} of the Olshanski $R$-matrix $S$.

\ble\label{lem:cross}
We have the relation
\beql{crosy}
(S^{-1})^{t_1}D_1S^{t_1}=D_1.
\eeq
\ele

\bpf
By regarding both sides of the relation as operators on the space $\CC^{n|n}\ot \CC^{n|n}$,
we will verify that the results of their applications to all basis vectors agree.
We have
\ben
S^{t_1}=\sum_{i,j=-n}^n q^{(\de_{i,j}+\de_{-i,j})(1-2\bj)}_{}E_{ii}\ot E_{jj}
+\ve \sum_{i>j}\big(E_{ji}(-1)^{\bi\bj}+E_{-j,-i}(-1)^{\bi\bj+\bi+\bj}\big)\ot E_{ji}.
\een
We find from \eqref{ssto} that $S^{-1}=S_{21}-\ve P$ and so
\ben
(S^{-1})^{t_1}=\sum_{i,j=-n}^n q^{(\de_{i,j}+\de_{-i,j})(1-2\bj)}_{}E_{jj}\ot E_{ii}
-\ve \sum_{i\leqslant j}E_{ij}\ot E_{ij}(-1)^{\bi\bj}
+\ve \sum_{i>j}E_{ij}\ot E_{-i,-j}(-1)^{\bi\bj}.
\een
These expressions show that the application of both sides of \eqref{crosy}
to a basis vector $e_i\ot e_j$ with $i\ne \pm j$ yields the same vector
$q^{-2\tss|i|}\tss e_i\ot e_j$.

Now assume that $j\in\{1,\dots,n\}$ and first apply both sides of \eqref{crosy}
to the basis vector $e_j\ot e_j$. The right hand side gives $q^{-2\tss j}\tss e_j\ot e_j$,
while for the left hand side we get
\ben
\bal
(S^{-1})^{t_1}D_1S^{t_1}(e_j\ot e_j)&
=(S^{-1})^{t_1}D_1\big(q\tss e_j\ot e_j+\ve\sum_{i<j}e_i\ot e_i\big)\\
&=(S^{-1})^{t_1}\big(q^{-2j+1}\tss e_j\ot e_j+\ve\sum_{i<j}q^{-2|i|}\tss e_i\ot e_i\big)
\eal
\een
which equals
\ben
q^{-2j+1}\big(q\tss e_j\ot e_j-\ve\sum_{i\leqslant j}e_i\ot e_i\big)+
\ve\sum_{i<j}q^{-2|i|}\big(q^{1-2\bi} e_i\ot e_i-\ve\sum_{k\leqslant i}(-1)^{\bi}\tss e_k\ot e_k\big).
\een
Here the coefficient of $e_j\ot e_j$ equals $q^{-2\tss j}$, whereas for $k<j$ the coefficient
of $e_k\ot e_k$ equals
\ben
-\ve\tss q^{-2j+1}+\ve\tss q^{1-2|k|-2\bk}-\ve^2\ts\sum_{i=k}^{j-1} q^{-2|i|}(-1)^{\bi}=0,
\een
thus completing the verification for the basis vector $e_j\ot e_j$.

As the next step, apply both sides of \eqref{crosy}
to the basis vector $e_j\ot e_{-j}$. The right hand side gives $q^{-2\tss j}\tss e_j\ot e_{-j}$,
while for the left hand side we get
\ben
\bal
(S^{-1})^{t_1}D_1S^{t_1}(e_j\ot e_{-j})&
=(S^{-1})^{t_1}D_1\big(q^{-1}\tss e_j\ot e_{-j}-\ve\sum_{i>j}e_i\ot e_{-i}\big)\\
&=(S^{-1})^{t_1}\big(q^{-2j-1}\tss e_j\ot e_{-j}-\ve\sum_{i>j}q^{-2\tss i}\tss e_i\ot e_{-i}\big)
\eal
\een
which equals
\ben
q^{-2j-1}\big(q\tss e_j\ot e_{-j}+\ve\sum_{i> j}e_i\ot e_{-i}\big)-
\ve\sum_{i>j}q^{-2\tss i}\big(q\tss e_i\ot e_{-i}+\ve\sum_{k> i}\tss e_k\ot e_{-k}\big).
\een
The coefficient of $e_j\ot e_{-j}$ equals $q^{-2\tss j}$, whereas for $k>j$ the coefficient
of $e_k\ot e_k$ equals
\ben
\ve\tss q^{-2j-1}-\ve\tss q^{-2k+1}-\ve^2\ts\sum_{i=j+1}^{k-1} q^{-2i}=0,
\een
thus completing the verification for the basis vector $e_j\ot e_{-j}$.

The remaining two cases of basis vectors $e_{-j}\ot e_j$ and $e_{-j}\ot e_{-j}$
are verified in the same way.
\epf

This completes the proof of relation \eqref{sutre} and the theorem.
\epf

\bre\label{rem:secre}\ \  (i) The same argument shows that
the elements
\beql{savagesrtce}
\str\ts D^{-1}(L^{-1}\oL)^k,\qquad k=1,2,\dots
\eeq
also belong to the center of the superalgebra $\U_q(\q_n)$. They coincide with
the central elements of \cite[Eg.~(8.7)]{s:qi} up to the factor $q^{2n+1}$.
It was pointed out in \cite[Corollary~3.4(b)]{s:qi} that all elements $(L_{11}\dots L_{nn})^k$
with $k\in\ZZ$
are central in $\U_q(\q_n)$ and that it is likely that together with the family \eqref{savagesrtce}
they generate the center of $\U_q(\q_n)$. One can also expect
that the elements \eqref{srtce} generate the center of the reflection equation
subalgebra $\U^{\circ}_q(\q_n)$; cf. a similar property for
the algebra $\U_q(\gl_n)$, as reviewed in \cite{jlm:qi}.

(ii) There is another crossing symmetry relation for the $R$-matrix $S$ which is proved
in the way similar to Lemma~\ref{lem:cross}:
\beql{crosytwo}
(S^{-1})^{t_2}D^{-1}_2S^{t_2}=D^{-1}_2,
\eeq
although we will not use it.

(iii) Note also that a general procedure to find a suitable matrix $D$ associated with
a {\em skew-invertible $R$-matrix} to satisfy the crossing symmetry relations
is described in \cite[Sec.~3.1]{i:qg}. However, it turned out to be possible
to verify the relations directly by taking $D$ in the form \eqref{d} by analogy with
the $R$-matrix associated with $\gl_n$; cf. \cite[Sec.~1]{rtf:ql}.
\qed
\ere

\section{Central elements in $\U_q(\wh\q_n)$}
\label{sec:cent-aff}

We will now work with the superalgebra $\U_q(\wh\q_n)$ as defined in Sec.~\ref{sec:eh}
and use the matrix $D$ defined in \eqref{d}.

\bth\label{thm:cnetaff}
There exists a series $z(u)$ in $u^{-1}$ with coefficients
in the superalgebra $\U_q(\wh\q_n)$ such that
\beql{lltr}
L(u)^tD(L(u)^{-1})^t=z(u)D.
\eeq
Moreover, the coefficients of the series $z(u)$ belong to the center of the superalgebra
$\U_q(\wh\q_n)$.
\eth

\bpf
We start by proving a crossing symmetry property for the Jones--Nazarov $R$-matrix,
analogous to Lemma~\ref{lem:cross}. Set $S_{21}(u,v)=PS(u,v)P$.

\ble\label{lem:croaff}
We have the relation
\beql{crosaffwa}
S(u,v)^{t_1}D^{-1}_1S_{21}(v,u)^{t_1}D_1=1
\eeq
in the superalgebra $\End\CC^{n|n}\ot\End\CC^{n|n}$.
\ele

\bpf
By \eqref{jns}
we have
\beql{suvt}
S(u,v)^{t_1}=S^{t_1}+\frac{\ve\tss P^{t_1}}{u^{-1}v-1}+\frac{\ve\tss (PJ_1J_2)^{t_1}}{u\tss v-1},
\eeq
while using \eqref{ssto} we can write
\beql{stovu}
S_{21}(v,u)^{t_1}=
(S^{-1})^{t_1}+\frac{\ve\tss v^{-1}u\tss P^{t_1}}{v^{-1}u-1}
+\frac{\ve\tss (J_1J_2 P)^{t_1}}{v\tss u-1}.
\eeq
Arguing as in the proof of Lemma~\ref{lem:cross}, apply the operators appearing
on both sides of
\eqref{crosaffwa} to all basis vectors $e_i\ot e_j$ of $\CC^{n|n}\ot\CC^{n|n}$.
We will use the expressions for $S^{t_1}$ and $(S^{-1})^{t_1}$ derived in that proof
along with the relation
\beql{crofin}
S^{t_1}D^{-1}_1(S^{-1})^{t_1}D_1=1
\eeq
which is equivalent to \eqref{crosy}. Note also the expressions
\ben
\bal
P^{t_1}&=\sum_{i,j=-n}^n E_{ji}\ot E_{ji}(-1)^{\bi\bj},\\
(PJ_1J_2)^{t_1}&=\sum_{i,j=-n}^n E_{-j,-i}\ot E_{ji}(-1)^{\bi\bj+\bi+\bj},\\
(J_1J_2P)^{t_1}&=\sum_{i,j=-n}^n E_{ji}\ot E_{-j,-i}(-1)^{\bi\bj},
\eal
\een
to be used in the calculations below.
Clearly, if $i\ne \pm j$, then the application of the operators
on both sides of \eqref{crosaffwa} to $e_i\ot e_j$ yields $e_i\ot e_j$ .
Now assume that $j\in\{1,\dots,n\}$ and apply the left hand side of \eqref{crosaffwa}
to the basis vector $e_j\ot e_j$.
Using \eqref{stovu} and \eqref{suvt}, expand the product to derive
from \eqref{ssto} that the result of the application is $e_j\ot e_j$
plus the sum of two expressions
\beql{phione}
S^{t_1}D^{-1}_1\Bigg(\frac{\ve\tss v^{-1}u\tss P^{t_1}}{v^{-1}u-1}
+\frac{\ve\tss (J_1J_2 P)^{t_1}}{v\tss u-1}\Bigg)D_1(e_j\ot e_j)
\eeq
and
\beql{phitwo}
\Bigg(\frac{\ve\tss P^{t_1}}{u^{-1}v-1}+\frac{\ve\tss (PJ_1J_2)^{t_1}}{u\tss v-1}\Bigg)
D^{-1}_1 S_{21}(v,u)^{t_1}D_1(e_j\ot e_j).
\eeq
Expanding \eqref{phione} we get
\begin{multline}
\frac{\ve\tss v^{-1}u}{v^{-1}u-1}\ts S^{t_1}D^{-1}_1P^{t_1}D_1(e_j\ot e_j)
=q^{-2\tss j}\ts \frac{\ve\tss v^{-1}u}{v^{-1}u-1}\ts S^{t_1}D^{-1}_1P^{t_1}(e_j\ot e_j)\\[0.4em]
=q^{-2\tss j}\ts \frac{\ve\tss v^{-1}u}{v^{-1}u-1}\ts S^{t_1}D^{-1}_1\sum_{i=-n}^n(e_i\ot e_i)
=q^{-2\tss j}\ts \frac{\ve\tss v^{-1}u}{v^{-1}u-1}\ts S^{t_1}\sum_{i=-n}^n\tss q^{2\tss|i|}(e_i\ot e_i)
\non
\end{multline}
which equals
\ben
q^{-2\tss j}\ts \frac{\ve\tss v^{-1}u}{v^{-1}u-1}\ts
\sum_{i=-n}^n\tss q^{2\tss|i|}\ts\Big(q^{1-2\tss \bi}\tss e_i\ot e_i+\ve\tss
\sum_{k<i}e_k\ot e_k\tss (-1)^{\bi}\Big)
\een
and simplifies to
\beql{sumone}
\frac{\ve\tss q^{2n-2j+1}\tss u}{u-v}\ts
\sum_{k=-n}^n e_k\ot e_k.
\eeq
Similarly, by expanding \eqref{phitwo} we obtain
\begin{multline}
\frac{\ve\tss P^{t_1}}{u^{-1}v-1}\ts D^{-1}_1\Bigg((S^{-1})^{t_1}
+\frac{\ve\tss v^{-1}u\tss P^{t_1}}{v^{-1}u-1}\Bigg)D_1(e_j\ot e_j)\\[0.4em]
=q^{-2\tss j}\ts\frac{\ve\tss P^{t_1}}{u^{-1}v-1}\ts D^{-1}_1
\Bigg(q^{-1}\tss e_j\ot e_j-\ve\tss \sum_{k<j}e_k\ot e_k
+\frac{\ve\tss v^{-1}u\tss }{v^{-1}u-1}\sum_{k=-n}^n e_k\ot e_k \Bigg).
\non
\end{multline}
Next, apply $D^{-1}_1$ and observe that for any $i\in\{-n,\dots,-1,1,\dots,n\}$,
\beql{ptone}
P^{t_1}(e_i\ot e_i)=(-1)^{\bi}\sum_{k=-n}^n e_k\ot e_k,
\eeq
to bring the expression to the form
\ben
q^{2\tss n-2\tss j+1}\ts\frac{\ve\tss P^{t_1}}{u^{-1}v-1}(e_j\ot e_j)
=\frac{\ve\tss q^{2n-2j+1}\tss u}{v-u}\ts
\sum_{k=-n}^n e_k\ot e_k
\een
so that its sum with the expression in \eqref{sumone} is zero thus completing
the verification for the vector $e_j\ot e_j$.

Now apply the operator on the left hand side of
\eqref{crosaffwa}
to the basis vector $e_j\ot e_{-j}$. As with the previous case, as a result, we get
$e_j\ot e_{-j}$ plus the sum of the expressions
\eqref{phione} and \eqref{phitwo}, where $e_j\ot e_j$ is replaced with $e_j\ot e_{-j}$.
This time, expanding \eqref{phione} we get
\begin{multline}
\frac{\ve}{vu-1}\ts S^{t_1}D^{-1}_1(J_1J_2P)^{t_1}D_1(e_j\ot e_{-j})
=q^{-2\tss j}\ts \frac{\ve}{vu-1}\ts S^{t_1}D^{-1}_1(J_1J_2P)^{t_1}(e_j\ot e_{-j})\\[0.4em]
=q^{-2\tss j}\ts \frac{\ve}{vu-1}\ts S^{t_1}D^{-1}_1\sum_{i=-n}^n(e_i\ot e_{-i})
=q^{-2\tss j}\ts \frac{\ve}{vu-1}\ts  S^{t_1}\sum_{i=-n}^n\tss q^{2\tss|i|}(e_i\ot e_{-i})
\non
\end{multline}
which equals
\ben
q^{-2\tss j}\ts \frac{\ve}{vu-1}\ts
\sum_{i=-n}^n\tss q^{2\tss|i|}\ts\Big(q^{-1+2\tss \bi}\tss e_i\ot e_{-i}-\ve\tss
\sum_{k>i}e_k\ot e_{-k}\tss (-1)^{\bi}\Big)
\een
and simplifies to
\beql{sumonepm}
\frac{\ve\tss q^{2n-2j+1}}{vu-1}\ts
\sum_{k=-n}^n e_k\ot e_{-k}.
\eeq
Working now with the respective expression in
\eqref{phitwo} we obtain
\begin{multline}
\frac{\ve\tss (PJ_1J_2)^{t_1}}{uv-1}\ts D^{-1}_1\Bigg((S^{-1})^{t_1}
+\frac{\ve\tss (J_1J_2P)^{t_1}}{vu-1}\Bigg)D_1(e_j\ot e_{-j})\\[0.4em]
=q^{-2\tss j}\ts\frac{\ve\tss (PJ_1J_2)^{t_1}}{uv-1}\ts D^{-1}_1
\Bigg(q\tss e_j\ot e_{-j}+\ve\tss \sum_{k>j}e_k\ot e_{-k}
+\frac{\ve}{vu-1}\sum_{k=-n}^n e_k\ot e_{-k} \Bigg).
\non
\end{multline}
Apply $D^{-1}_1$ and observe that for any $i\in\{-n,\dots,-1,1,\dots,n\}$,
\ben
(PJ_1J_2)^{t_1}(e_i\ot e_{-i})=(-1)^{\bi+1}\sum_{k=-n}^n e_k\ot e_{-k},
\een
to bring the expression to the form
\ben
q^{2\tss n-2\tss j+1}\ts\frac{\ve\tss (PJ_1J_2)^{t_1}}{uv-1}(e_j\ot e_{-j})
=-\frac{\ve\tss q^{2n-2j+1}}{uv-1}\ts
\sum_{k=-n}^n e_k\ot e_{-k}
\een
so that its sum with the expression in \eqref{sumonepm} is zero thus completing
the verification for the vector $e_j\ot e_{-j}$.

The calculations in the remaining two cases of basis vectors $e_{-j}\ot e_j$ and $e_{-j}\ot e_{-j}$
are quite similar and will be omitted.
\epf

Returning to the proof of the theorem, multiply both sides of \eqref{RTTuv}
by $L_1(u)^{-1}$ from the left and from the right, then apply the
super-transposition $t_1$ to get
\ben
L_2(v)\tss (L_1(u)^{-1})^{t_1}S(u,v)^{t_1}=S(u,v)^{t_1}\tss (L_1(u)^{-1})^{t_1}L_2(v).
\een
This implies
\beql{stll}
(S(u,v)^{t_1})^{-1}\tss L_2(v)\tss (L_1(u)^{-1})^{t_1}
=(L_1(u)^{-1})^{t_1}L_2(v)(S(u,v)^{t_1})^{-1}.
\eeq
Now apply Lemma~\ref{lem:croaff}
to write
\ben
(S(u,v)^{t_1})^{-1}=D^{-1}_1S_{21}(v,u)^{t_1}D_1.
\een
Hence, we derive from \eqref{stll} that
\ben
D^{-1}_1S_{21}(v,u)^{t_1}D_1\tss L_2(v)\tss (L_1(u)^{-1})^{t_1}=
(L_1(u)^{-1})^{t_1}L_2(v)\tss D^{-1}_1S_{21}(v,u)^{t_1}D_1.
\een
Due to \eqref{stovu},
multiplying both sides by $v^{-1}u-1$ and setting $v=u$ we get
\ben
D^{-1}_1P^{t_1}D_1\tss L_2(u)\tss (L_1(u)^{-1})^{t_1}=
(L_1(u)^{-1})^{t_1}L_2(u)\tss D^{-1}_1P^{t_1}D_1
\een
which is equivalent to
\beql{ptdll}
P^{t_1}D_1\tss L_2(u)\tss (L_1(u)^{-1})^{t_1}D^{-1}_1=
D_1(L_1(u)^{-1})^{t_1}L_2(u)\tss D^{-1}_1P^{t_1}.
\eeq
As relation \eqref{ptone} shows, the image of the operator
$P^{t_1}$ in $\CC^{n|n}\ot\CC^{n|n}$ is one-dimensional so that both sides
of \eqref{ptdll} must equal to $P^{t_1}$ times a series $z(u)$ in $u^{-1}$
with coefficients in the superalgebra $\U_q(\wh\q_n)$. Since
$P^{t_1}L_2(u)=P^{t_1}L_1(u)^{t_1}$, we can write
\ben
P^{t_1}L_1(u)^{t_1} \tss D_1\tss (L_1(u)^{-1})^{t_1}D^{-1}_1=P^{t_1}\tss z(u).
\een
Note that $\str_2\tss P^{t_1}=1$ so that
by taking the supertrace over the second copy of $\End\CC^{n|n}$, we get
\ben
L_1(u)^{t_1} \tss D_1\tss (L_1(u)^{-1})^{t_1}D^{-1}_1=z(u),
\een
which is equivalent to \eqref{lltr}.

Finally, we will show that all coefficients of $z(u)$ are central in $\U_q(\wh\q_n)$.
Working in the superalgebra $\U_q(\wh\q_n)\ot \End\CC^{n|n}\ot \End\CC^{n|n}$
write
\ben
z(u)\tss L_2(v)=L_1(u)^{t_1} \tss D_1\tss (L_1(u)^{-1})^{t_1}D^{-1}_1 \tss L_2(v)
=L_1(u)^{t_1} \tss D_1\tss (L_1(u)^{-1})^{t_1} \tss L_2(v)D^{-1}_1.
\een
Now apply \eqref{stll} to write this expression as
\beql{lds}
L_1(u)^{t_1} \tss D_1\tss (S(u,v)^{t_1})^{-1}\tss L_2(v)(L_1(u)^{-1})^{t_1}S(u,v)^{t_1}D^{-1}_1.
\eeq
By Lemma~\ref{lem:croaff},
\ben
D_1(S(u,v)^{t_1})^{-1}=S_{21}(v,u)^{t_1}D_1.
\een
On the other hand, \cite[Lemma~4.1]{jn:as} implies that
\beql{jnlem}
S_{21}(v,u)=S(u,v)^{-1}\ts A(u,v),
\eeq
where
\ben
A(u,v)=1-\frac{\ve^2\ts uv}{(u-v)^2}-\frac{\ve^2\ts uv}{(uv-1)^2}.
\een
Therefore, expression \eqref{lds} equals
\beql{alds}
A(u,v)\ts L_1(u)^{t_1} \tss  (S(u,v)^{-1})^{t_1}\tss L_2(v)D_1(L_1(u)^{-1})^{t_1}S(u,v)^{t_1}D^{-1}_1.
\eeq
We have the relation
\ben
L_1(u)^{t_1} \tss  (S(u,v)^{-1})^{t_1}\tss L_2(v)=L_2(v) \tss  (S(u,v)^{-1})^{t_1}\tss L_1(u)^{t_1},
\een
which follows from \eqref{RTTuv} by multiplying both sides
by $S(u,v)^{-1}$ from the left and from the right and applying the
super-transposition $t_1$. Hence \eqref{alds} can be written as
\beql{aldslsv}
A(u,v)\tss L_2(v) \tss  (S(u,v)^{-1})^{t_1}\ts L_1(u)^{t_1}D_1(L_1(u)^{-1})^{t_1}S(u,v)^{t_1}D^{-1}_1
\eeq
which equals
\ben
A(u,v)\tss L_2(v) \tss  (S(u,v)^{-1})^{t_1}\ts D_1\ts z(u)\ts S(u,v)^{t_1}D^{-1}_1
=L_2(v) \tss z(u)\tss S_{21}(v,u)^{t_1}\ts D_1\ts S(u,v)^{t_1}D^{-1}_1
\een
by \eqref{lltr} and \eqref{jnlem}. The last expression equals $L_2(v) \tss z(u)$
by Lemma~\ref{lem:croaff}, thus completing the proof of the theorem.
\epf

\bre\label{rem:secreaff}
Another crossing symmetry relation for the $R$-matrix $S(u,v)$ takes the form
\ben
S(u,v)^{t_2}D_2S_{21}(v,u)^{t_2}D^{-1}_2=1
\een
which is verified
in a way similar to Lemma~\ref{lem:croaff}. Note also that \eqref{lltr}
has an equivalent form
\ben
(L(u)^{-1})^t D^{-1}L(u)^t=z(u)D^{-1}
\een
which follows from \eqref{lltr} or \eqref{ptdll}.
\qed
\ere

For a concluding remark, note that one can expect that the images of
the coefficients of the series $z(u)$ under the evaluation homomorphism
of Theorem~\ref{thm:eval} are related to the central elements
of $\U_q(\q_n)$ constructed in Theorem~\ref{thm:centfin}.
This is supported by the simplest case of $n=1$ considered below, but
we do not have a proof for general values of $n$.

\bex\label{ex:none}
In the superalgebra $\U_q(\q_1)$ we have the relations $L_{11}L_{-1,1}=L_{-1,1}L_{11}$
so that the elements $L_{11}$ and $L_{-1,-1}=L^{-1}_{11}$ lie in the center of $\U_q(\q_1)$.
Furthermore, we have
\ben
L_{-1,1}^2=\frac{q^2-1}{q^2+1}\ts \big(L_{11}^{2}-L_{-1,-1}^2\big).
\een
The matrices $L$ and $\oL$ have the form
\ben
L=\begin{bmatrix}L^{}_{-1,-1}&L^{}_{-1,1}\\
0& L^{}_{11}
\end{bmatrix}
\Fand
\oL=\begin{bmatrix}-L^{}_{11}&0\\
\phantom{-}L^{}_{-1,1}& -L^{}_{-1,-1}
\end{bmatrix}.
\een
Hence,
\ben
M=L\tss \oL^{-1}=
\begin{bmatrix}-L^{2}_{-1,-1}-L^{2}_{-1,1}&-L^{}_{-1,1}L_{11}\\[0.3em]
-L^{}_{-1,1}L_{11}& -L^{2}_{11}
\end{bmatrix}.
\een
Up to scalar factors,
the central elements of $\U_q(\q_1)$ provided by Theorem~\ref{thm:centfin}
are found as the supertraces of the powers of $M$.
By calculating the powers,
for all $k\geqslant 0$ we get
\ben
\str\ts M^{k+1}=(-1)^{k}\big(L^{2}_{11}-L^{2}_{-1,-1}-L^{2}_{-1,1}\big)
\sum_{r\geqslant 0}(-1)^r\binom{k-r}{r}\big(L^{2}_{11}+L^{2}_{-1,-1}+L^{2}_{-1,1}\big)^{k-2r}.
\een

On the other hand,
the image of the series $z(u)$ under the evaluation homomorphism provided by Theorem~\ref{thm:eval}
takes the form
\ben
\ev:z(u)\mapsto \big(L+\oL\tss u^{-1}\big)^t\ts \Big(\big(L+\oL\tss u^{-1}\big)^{-1}\Big)^t.
\een
By taking the $(1,1)$ entry we get
\ben
\ev:z(u)\mapsto \Big(1-\big(L^{2}_{11}+L^{2}_{-1,-1}-L^{2}_{-1,1}\big)\ts u^{-1}+u^{-2}\Big)
\Big(1-\big(L^{2}_{11}+L^{2}_{-1,-1}+L^{2}_{-1,1}\big)\ts u^{-1}+u^{-2}\Big)^{-1}.
\een
A straightforward calculation shows that
\ben
\ev:z(u)\mapsto 1+(1-q^2)\ts \sum_{k=1}^{\infty}(-1)^{k}\ts \str\ts M^k\tss u^{-k}.
\een
\eex

\section*{Declarations}

\subsection*{Competing interests}
The authors have no competing interests to declare that are relevant to the content of this article.

\subsection*{Acknowledgements}
This project was initiated during the second named author's visit
to the Central China Normal University in Wuhan. He
is grateful to the School of Mathematics and Statistics for warm hospitality.
His work was also supported by the Australian Research Council, grant DP240101572.
Liu is supported by the National Natural Science Foundation of China (Grant No. 12471026).
Zhang is supported by the National Natural Science Foundation of China (Grant No. 12571026).

\subsection*{Availability of data and materials}
No data was used for the research described in the article.


\bigskip
\bigskip

\small

\noindent
M.L.:\newline
School of Artificial Intelligence\newline
Jianghan University,
Wuhan 430056, Hubei, China\newline
{\tt ming.l1984@gmail.com}

\vspace{5 mm}

\noindent
A.M.:\newline
School of Mathematics and Statistics\newline
University of Sydney,
NSW 2006, Australia\newline
{\tt alexander.molev@sydney.edu.au}

\vspace{5 mm}

\noindent
J.Z.:\newline
School of Mathematics and Statistics and Hubei Key Lab - Math. Sci.\newline
Central China Normal University, Wuhan 430079, Hubei, China\newline
{\tt jzhang@ccnu.edu.cn}

\end{document}